\documentclass{birkmult}
\usepackage{amssymb}

\usepackage{mathrsfs}
\usepackage{multicol}
\usepackage{comment}

 \newtheorem{thm}{Theorem}[section]
 
 \newtheorem{lem}[thm]{Lemma}
 
 \theoremstyle{definition}
 \newtheorem{defn}[thm]{Definition}
 \theoremstyle{remark}
 \newtheorem{rem}[thm]{Remark}
 
  \newtheorem{cla}[thm]{Claim}
 \numberwithin{equation}{section}

\begin{document}
%
%
%
%
%
%
%
%
%
\title[Calabi
product Lagrangian immersions in $\mathbb{CP}^n$ and
$\mathbb{CH}^n$]{\large
 Calabi
product Lagrangian immersions in complex projective space and
complex hyperbolic space }
\author[Haizhong Li]{Haizhong Li}

\address{%
Department of Mathematical Sciences\\
Tsinghua University\\
100084, Beijing\\
People's Republic of China}

\email{hli@math.tsinghua.edu.cn}

\thanks{The first author was supported by NSFC grant No. 10971110 and Tsinghua University--K.U.Leuven Bilateral
scientific cooperation Fund.
          The second author was supported by NSFC grant No. 10701007 and
Tsinghua University--K.U.Leuven Bilateral scientific cooperation
Fund.}
\author{Xianfeng Wang}
\address{Department of Mathematical Sciences\br
Tsinghua University\br 100084, Beijing\br People's Republic of
China} \email{xf-wang06@mails.tsinghua.edu.cn}


\date{September 16, 2010}
\dedicatory{Dedicated to Heinrich Wefelscheid on the occasion of his
70th birthday}


\begin{abstract}
Starting from two Lagrangian immersions and a Legendre curve
$\tilde{\gamma}(t)$ in $\mathbb{S}^3(1)$ (or in
$\mathbb{H}_1^3(1)$), it is possible to construct a new Lagrangian
immersion in $\mathbb{CP}^n$ (or in $\mathbb{CH}^n$), which is
called a warped product Lagrangian immersion. When
$\tilde{\gamma}(t)=(r_1e^{i(\frac{r_2}{r_1}at)}, r_2e^{i(-
\frac{r_1}{r_2}at)})$ (or
$\tilde{\gamma}(t)=(r_1e^{i(\frac{r_2}{r_1}at)}, r_2e^{i(
\frac{r_1}{r_2}at)})$), where $r_1$, $r_2$, and $a$ are positive
constants with $r_1^2+r_2^2=1$ (or $-r_1^2+r_2^2=-1$), we call the
new Lagrangian immersion a Calabi product Lagrangian immersion. In
this paper, we study the inverse problem: how to determine from the
properties of the second fundamental form whether a given Lagrangian
immersion of $\mathbb{CP}^n$ or $\mathbb{CH}^n$ is a Calabi product
Lagrangian immersion. When the Calabi product is minimal, or is
Hamiltonian minimal, or has parallel second fundamental form, we
give some further characterizations.
\end{abstract}

\maketitle
\section{Introduction}

Let $\psi:M\to \bar{M}^n$ be an isometric immersion from an
n-dimensional Riemannian manifold into a complex n-dimensional
K\"{a}hler manifold $\bar{M}^n$. $M$ is called a
\textit{Lagrangian submanifold} if the almost complex structure $J$
of $\bar{M}^n$ carries each tangent space of $M$ into its
corresponding normal space.

Let $\bar{M}^n(4c)$ denote the complex space form with constant
holomorphic sectional curvature $4c$.

When $c>0$, $\bar{M}^n(4c)=\mathbb{CP}^n(4c)$.
Let
$\mathbb{S}^{2n+1}(c)=\{z\in \mathbb{C}^{n+1}: \langle z,z \rangle =\frac{1}{c}>0\}$
be the hypersphere of $\mathbb{C}^{n+1}$ with constant sectional curvature $c$ centered at the origin.
We have the Hopf fibration $\Pi:\mathbb{S}^{2n+1}(c)\to \mathbb{CP}^n(4c)$ (see \cite{R}).
On $\mathbb{S}^{2n+1}(c)$ we consider the contact structure $\phi$ (i.e., the projection of the complex
structure $J$ of $\mathbb{C}^{n+1}$ on the tangent bundle of $\mathbb{S}^{2n+1}(c)$) and the structure vector field
$\zeta=Jx$, where $x$ is the position vector.
An isometric immersion $f: M\to \mathbb{S}^{2n+1}(c)$ is
called\textit{ C-totally real} if $\zeta$ is normal to $f_{*}(TM)$ and $\langle \phi(f_{*}(TM)),f_{*}(TM)\rangle=0$, where
$\langle,\rangle$ denotes the inner product on $\mathbb{C}^{n+1}$.

When $c<0$, $\bar{M}^n(4c)=\mathbb{CH}^n(4c)$. We consider the
$(n+1)$-dimensional complex number space $\mathbb{C}^{n+1}_1$
endowed with the pseudo-Euclidean metric $g_1$ given by
$g_1=-dz_1d\bar{z}_1+\sum\limits_{j=2}^{n+1}dz_jd\bar{z}_j.$ Let
$\mathbb{H}_1^{2n+1}(c)=\{z\in\mathbb{C}^{n+1}_1: \langle z,z
\rangle _1=\frac{1}{c}<0\}$ be the anti-de Sitter space-time, where
$\langle,\rangle_1$ denotes the inner product on
$\mathbb{C}_1^{n+1}$. We have the Hopf fibration
$\Pi:\mathbb{H}_1^{2n+1}(c)\to \mathbb{CH}^n(4c).$ On
$\mathbb{H}_1^{2n+1}(c)$ we consider the contact structure $\phi$
(i.e., the projection of the complex structure $J$ of
$\mathbb{C}_1^{n+1}$ on the tangent bundle of
$\mathbb{H}_1^{2n+1}(c)$) and the structure vector field $\zeta=Jx$,
where $x$ is the position vector. An isometric immersion $f: M\to
\mathbb{H}_1^{2n+1}(c)$ is called\textit{ C-totally real} if $\zeta$
is normal to $f_{*}(TM)$ and $\langle
\phi(f_{*}(TM)),f_{*}(TM)\rangle_1=0$.

C-totally real curves in $\mathbb{S}^3(c)$ (or in $\mathbb {H}_1^3(c)$) are called \textit{Legendre curves}.
It is known that $\tilde{\gamma}(t):I\to \mathbb{S}^3(c)$ (or $\tilde{\gamma}(t):I\to \mathbb{H}_1^3(c)$ is a
Legendre curve if and only if $\tilde{\gamma}(t)$ satisfies the horizontality condition: $\langle \tilde{\gamma}',J\tilde{\gamma}\rangle=0$ (or
$\langle \tilde{\gamma}',J\tilde{\gamma}\rangle_1=0$) (see \cite{Ch1}).

In \cite{BMV} and \cite{MV}, the authors introduced a method to construct a new Lagrangian immersion from two lower dimensional
Lagrangian immersions in complex projective space, they called the new immersion a warped product Lagrangian immersion.
This construction slightly
generalized the construction in \cite{CLU} which is
analogous to the well-known Calabi product in affine differential
geometry (see \cite{HLV}). They gave some characterizations of the warped product Lagrangian immersion
constructed from two lower dimensional Lagrangian immersions.

\begin{defn}\label{def1.1}
Let $\psi_i:(M_i,g_i)\to \mathbb{CP}^{n_i}(4)$, $i=1,2,$ be two Lagrangian immersions and
let $\tilde{\gamma}=(\tilde{\gamma}_1,\tilde{\gamma}_2): I\to \mathbb{S}^3(1)\subset
\mathbb{C}^2$ be a Legendre curve. Then
$\psi=\Pi(\tilde{\gamma}_1\tilde{\psi}_1;\tilde{\gamma}_2\tilde{\psi}_2):I\times
M_1\times M_2\to \mathbb{CP}^{n}(4)$ is a Lagrangian immersion, where $n=n_1+n_2+1$,
$\tilde{\psi}_i:M_i\to \mathbb{S}^{2n_i+1}(1)$ are horizontal lifts of $\psi_i$,
$i=1,2$, respectively and $\Pi$ is the Hopf fibration. We call $\psi$ a \textit{warped
product} Lagrangian immersion of $\psi_1$ and $\psi_2$. When $n_1$ (or $n_2$) is zero, we
call $\psi$ a \textit{warped product} Lagrangian immersion of $\psi_2$ (or $\psi_1$)  and
a point.
\end{defn}
\begin{defn}\label{def1.2}
In Definition \ref{def1.1}, when
$\tilde{\gamma}(t)=(r_1e^{i(\frac{r_2}{r_1}at)}, r_2e^{i(-
\frac{r_1}{r_2}at)})$, where $r_1$, $r_2$, and $a$ are positive
constants with $r_1^2+r_2^2=1$,
 we call $\psi$ a
\textit{Calabi product} Lagrangian immersion of $\psi_1$ and
$\psi_2$. When $n_1$ (or $n_2$) is zero, we call  $\psi$ a
\textit{Calabi product} Lagrangian immersion of $\psi_2$ (or
$\psi_1$)  and a point.
\end{defn}

Analogously, we define warped product Lagrangian immersion and Calabi
product Lagrangian immersion in complex hyperbolic space.

\begin{defn}\label{def1.4}
Let $\psi_1:(M_1,g_1)\to \mathbb{CH}^{n_1}(-4)$ and  $\psi_2:(M_2,g_2)\to
\mathbb{CP}^{n_2}(4)$ be two Lagrangian immersions, let $\tilde{\psi}_1:M_1\to
\mathbb{H}_1^{2n_1+1}(-1)$ and $\tilde{\psi}_2:M_2\to \mathbb{S}^{2n_2+1}(1)$ be
horizontal lifts of $\psi_1$ and $\psi_2$ respectively and let
$\tilde{\gamma}=(\tilde{\gamma}_1,\tilde{\gamma}_2): I\to \mathbb{H}_1^3(1)\subset
\mathbb{C}_1^2$ be a Legendre curve. Then
$\psi=\Pi(\tilde{\gamma}_1\tilde{\psi}_1;\tilde{\gamma}_2\tilde{\psi}_2):I\times
M_1\times M_2\to \mathbb{CH}^{n}(-4)$ is a Lagrangian immersion, where $n=n_1+n_2+1$ and
$\Pi$ is the Hopf fibration. We call $\psi$ a \textit{warped product} Lagrangian
immersion of $\psi_1$ and $\psi_2$. When $n_1$ (or $n_2$) is zero, we call  $\psi$ a
\textit{warped product} Lagrangian immersion of $\psi_2$ (or $\psi_1$)  and a point.
\end{defn}
\begin{defn}\label{def1.3}
In Definition \ref{def1.4}, when
$\tilde{\gamma}(t)=(r_1e^{i(\frac{r_2}{r_1}at)}, r_2e^{i(
\frac{r_1}{r_2}at)})$, where $r_1$, $r_2$, and $a$ are positive
constants with  $-r_1^2+r_2^2=-1$, we call $\psi$ a \textit{Calabi
product} Lagrangian immersion of $\psi_1$ and $\psi_2$. When $n_1$
(or $n_2$) is zero, we call  $\psi$ a \textit{Calabi product}
Lagrangian immersion of $\psi_2$ (or $\psi_1$)  and a point.
\end{defn}

In this paper, we study the inverse problem: how to
determine from the properties of the second fundamental form whether
a given Lagrangian immersion of $\mathbb{CP}^n$ or $\mathbb{CH}^n$
is a Calabi product Lagrangian immersion. When the Calabi product is
minimal, or is Hamiltonian minimal, or has parallel second fundamental
form, we give some further characterizations. We have the following
results:

\begin{thm}\label{thm1.3}
Let $\psi:M\to \mathbb{CP}^n(4)$ be a Lagrangian immersion, $\psi$
is locally a Calabi product Lagrangian immersion of an
(n-1)-dimensional Lagrangian immersion $\psi_1: M_1\to
\mathbb{CP}^{n-1}(4)$ and a point if and only if $M$ admits two
orthogonal distributions $\mathcal{D}_1$ (of dimension 1, spanned by
a unit vector $E_1$) and $\mathcal{D}_2$ (of dimension $n-1$,
spanned by $\{E_2,\ldots, E_n\}$), and there exist two real
constants $\lambda_1$ and $\lambda_2$ such that
\begin{equation}
\left\{
\begin{aligned}
&h(E_1,E_1)=\lambda_1JE_1,~h(E_1,E_i)=\lambda_2JE_i,\\
&\lambda_1\neq 2 \lambda_2,~i=2,\ldots,n.
\end{aligned}
\right.
\label{1.1}
\end{equation}
Moreover, a Lagrangian immersion  $\psi:M\to \mathbb{CP}^n(4)$ satisfying the above conditions has the following properties:

{\rm(i)} $\psi$ is  Hamiltonian minimal if and only if $\psi_1$ is
Hamiltonian minimal.

{\rm(i)} $\psi$ is minimal if and only if $\lambda_2=\pm \frac{1}{\sqrt{n}}$ and $\psi_1$ is minimal. In this case,
up to a reparametrization and a rigid motion of $\mathbb{CP}^n$, locally we have $M=I\times M_1$ and
$\psi$ is given by
$\psi=\Pi\circ\tilde{\psi}$ with
$$
\tilde{\psi}(t,p)=(\sqrt{\frac{n}{n+1}}e^{i\frac{1}{n+1}t}\tilde{\psi}_1(p),
\sqrt{\frac{1}{n+1}}e^{-i\frac{n}{n+1}t}),~(t,p)\in I\times M_1,
$$
where $\Pi$ is the Hopf fibration and $\tilde{\psi}_1:M_1\to
\mathbb{S}^{2n-1}(1)$ is the horizontal lift of $\psi_1$.
\end{thm}

\begin{thm}\label{thm1.4}
Let $\psi:M\to \mathbb{CP}^n(4)$ be a Lagrangian immersion. If $M$
admits two orthogonal distributions $\mathcal{D}_1$ (of dimension 1,
spanned by a unit vector $E_1$) and $\mathcal{D}_2$ (of dimension
$n-1$, spanned by $\{E_2,\ldots, E_n\}$), and there exist local
functions $\lambda_1,~\lambda_2$ such that \eqref{1.1} holds, then
$M$ has parallel second fundamental form if and only if $\psi$ is
locally a Calabi product Lagrangian immersion of a point and an
$(n-1)$-dimensional Lagrangian immersion $\psi_1: M_1\to
\mathbb{CP}^{n-1}(4)$ which has parallel second fundamental form.
\end{thm}

\begin{thm}\label{thm1.5}
Let $\psi:M\to \mathbb{CH}^n(-4)$ be a Lagrangian immersion, $\psi$
is locally a Calabi product Lagrangian immersion of an
(n-1)-dimensional Lagrangian immersion $\psi_1: M_1\to
\mathbb{CP}^{n-1}(4)$ (or $\psi_2: M_2\to \mathbb{CH}^{n-1}(-4)$)
and a point if and only if
 $M$ admits two orthogonal
distributions $\mathcal{D}_1$ (of dimension 1, spanned by a unit vector $E_1$) and
$\mathcal{D}_2$ (of dimension $n-1$, spanned by $\{E_2,\ldots, E_n\}$), and there exist
two real constants $\lambda_1$ and $\lambda_2$  such that \eqref{1.1} holds. Moreover, a
Lagrangian immersion  $\psi:M\to \mathbb{CH}^n(-4)$ satisfying the above conditions can
not be minimal and $\psi$ is Hamiltonian minimal if and only if $\psi_1$ (or $\psi_2$) is
Hamiltonian minimal.
\end{thm}

 The paper is organized as follows. In Section 2 we recall the basic formulas for Lagrangian submanifolds of
complex space forms. In Section 3 we give some properties of the Lagrangian immersions in
the complex space forms with the second fundamental form satisfying \eqref{1.1}. In
Section 4 and Section 5, we
 give some characterizations of the warped product Lagrangian immersions
and the Calabi product Lagrangian immersions in complex projective space and
complex hyperbolic space respectively.

\section{Preliminaries}
In this section, $M$ will always denote an n-dimensional
Lagrangian submanifold of $\bar{M}^n(4c)$ which is an n-dimensional
complex space form with constant holomorphic sectional curvature
$4c$. We denote the Levi-Civita connections on $M$,
$\bar{M}^n(4c)$ and the normal bundle by $\nabla$, $D$ and
$\nabla_X^{\bot}$ respectively. The formulas of Gauss and Weingarten
are given by (see  \cite{CO})
$$
D_XY=\nabla_XY+h(X,Y), ~D_X\xi=-A_{\xi}X+\nabla_X^{\bot}\xi,
$$
where $X$ and $Y$ are tangent vector fields and $\xi$ is a normal
vector field on $M$.

The Lagrangian condition implies that (see \cite{Ch1})
$$\nabla_X^{\bot}JY=J\nabla_XY, ~A_{JX}Y=-J h(X,Y)=A_{JY}X,$$
where $h$ is the second fundamental form and $A$ denotes the shape
operator.

We denote the curvature tensors of $\nabla$ and $\nabla_X^{\bot}$ by
$R$ and $R^{\bot}$ respectively. The first covariant
derivative of $h$ is defined by
\begin{equation}
\begin{aligned}
(\nabla h)(X,Y,Z)=\nabla^{\bot}_X h(Y,Z)-h(\nabla_X
Y,Z)-h(\nabla_X Z,Y),\label{23}
\end{aligned}
\end{equation}
where $X$, $Y$, $Z$ and $W$ are tangent vector fields.

The equations of Gauss, Codazzi and Ricci for a Lagrangian
submanifold of $\bar{M}^n(4c)$ are given by (see \cite{CO})

\begin{equation}\begin{aligned} \langle R(X,Y)Z,&W \rangle ~= \langle h(Y,Z),h(X,W) \rangle - \langle h(X,Z),h(Y,W) \rangle \\
&+c~( \langle X,W \rangle  \langle Y,Z \rangle - \langle X,Z \rangle  \langle Y,W \rangle ),\label{24}\end{aligned}\end{equation}
\begin{equation}(\nabla h)(X,Y,Z)=(\nabla h)(Y,X,Z),\label{25}\end{equation}
\begin{equation}\begin{aligned} \langle &R^{\bot}(X,Y)JZ,JW \rangle ~= \langle [A_{JZ},A_{JW}]X,Y \rangle \\
&+c~( \langle X,W \rangle  \langle Y,Z \rangle - \langle X,Z \rangle  \langle Y,W \rangle ),\label{26}\end{aligned}\end{equation}
where $X$, $Y$ $Z$ and $W$ are tangent vector fields. Note that for
a Lagrangian submanifold the equations of Gauss and Ricci are
mutually equivalent.

The Lagrangian condition implies that
\begin{equation} \langle R^{\bot}(X,Y)JZ,JW \rangle ~=~ \langle R(X,Y)Z,W \rangle ,\label{28}\end{equation}
\begin{equation} \langle h(X,Y),JZ \rangle ~=~ \langle h(X,Z),JY \rangle ,\label{29}\end{equation}
for tangent vector fields $X$, $Y$, $Z$ and $W$.
From \eqref{23} and \eqref{29}, we also have
\begin{equation} \langle (\nabla h)(W,X,Y),JZ \rangle ~=~ \langle (\nabla h)(W,X,Z),JY \rangle ,\label{295}\end{equation}
for tangent vector fields $X$, $Y$, $Z$ and $W$.

\section{Some lemmas}
We assume that $M$ is a Lagrangian submanifold in $\bar{M}^n(4c)$,
$M$ admits two orthogonal distributions $\mathcal{D}_1$ (of
dimension 1, spanned by a unit vector $E_1$) and $\mathcal{D}_2$ (of
dimension $n-1$, spanned by $\{E_2,\ldots, E_n\}$), and there exist
local functions $\lambda_1,~\lambda_2$ such that \eqref{1.1} holds.
In the following we always assume that $2\leq i,j,k\leq n$.

From Lemma 2.1 and Lemma 2.4 of \cite{BMV}, we have the following
lemmas:
\begin{lem}\label{lemma3.1}We have
$$
\left\{
\begin{aligned}
&E_i(\lambda_1)=(\lambda_1-2\lambda_2) \langle \nabla_{E_1}E_1,E_i \rangle ,\\
&(\lambda_1-2\lambda_2) \langle \nabla_{E_i}E_j,E_1 \rangle +E_1(\lambda_2)\delta_{ij}- \langle h(E_i,E_j),J\nabla_{E_1}E_1 \rangle =0.
\end{aligned}
\right.
$$
\end{lem}

\begin{lem}\label{lemma3.2}We have
$$
\begin{aligned}
&E_i(\lambda_2)\delta_{jk}- \langle h(E_j.E_k),J\nabla_{E_i}E_1 \rangle \\
&=E_1( \langle h(E_i,E_j),JE_k \rangle )-\sigma( \langle h(E_i,E_j),J\nabla_{E_1}E_k \rangle ).
\end{aligned}
$$
\end{lem}

From now on, we assume that $\nabla$$\lambda_1$ lies in $\mathcal{D}_1$ and $n\geq3$. From \eqref{1.1} we know that
$\lambda_1-2\lambda_2\neq 0$, so from the first equation of Lemma \ref{lemma3.1} we have
$
\nabla_{E_1}E_1=0.
$
which together with the second equation of Lemma \ref{lemma3.1} imply
\begin{equation}
 \langle \nabla_{E_i}E_j,E_1 \rangle =-\frac{E_1(\lambda_2)}{\lambda_1-2\lambda_2}\delta_{ij}.
\label{3.9}
\end{equation}

Since $
\nabla_{E_1}E_1=0,
$ we have $\mathcal{D}_1$ is integrable and the integral curves of $E_1$ are geodesics of $M$.
By \eqref{3.9} we have $ \langle [E_i,E_j],E_1 \rangle =0$, which implies $D_2$
is integrable. It is sufficient to conclude that locally $M$ is
isometric with $I\times M_1$ where $D_1$ is tangent to $I$ and $D_2$
is tangent to $M_1$. The product structure of $M$ implies the
existence of local coordinate $(t,p)$ for $M$ based on an open
subset containing the origin of $R\times R^{n_1}$, such that $D_1$
is given by $dp=0$ with $E_1=\frac{\partial}{\partial t}$ and $D_2$ is given by $dt=0$.

Let $k=\frac{E_1(\lambda_2)}{\lambda_1-2\lambda_2}$,
\eqref{3.9} implies $\nabla_{E_i}E_1=kE_i$, which together with Lemma \ref{lemma3.2} imply that
$
E_i(\lambda_2)\delta_{jk}=E_j(\lambda_2)\delta_{ik}.$
Since $n\geq 3$, we have $E_i(\lambda_2)=0,$ which implies that $\nabla$$\lambda_2$ and $\nabla$$k$ also lie in $\mathcal{D}_1$.
Hence we have $\lambda_1,~\lambda_2$ and $k$ are
 all functions of $t$.

 Let  $\lambda_1',~\lambda_2'$ and $k'$ denote the derivatives of $\lambda_1,~\lambda_2$ and $k$ with respect to $t$
 respectively. By the definition of the curvature tensor, we have
$$
\begin{aligned}
&\langle R(E_1,E_i)E_1,E_i\rangle= \langle \nabla_{E_1}\nabla_{E_i}E_1,E_i \rangle - \langle \nabla_{E_i}\nabla_{E_1}E_1,E_i \rangle \\
&- \langle \nabla_{\nabla_{E_1}E_i}E_1,E_i \rangle
+ \langle \nabla_{\nabla_{E_i}E_1}E_1,E_i \rangle=k'+k^2.
\end{aligned}
$$
On the other hand, by Gauss equation \eqref{24} we have
$$
\begin{aligned}
\langle R(E_1,E_i)E_1,E_i\rangle&=-c+ \langle h(E_i,E_1),h(E_i,E_1) \rangle - \langle h(E_i,E_i),h(E_1,E_1) \rangle \\
&=-c+\lambda_2^2-\lambda_1\lambda_2.
\end{aligned}
$$

By comparing both expressions, we have (c.f \cite{Ch1})
\begin{equation}
k'+k^2+\lambda_1\lambda_2-\lambda_2^2+c=0.\label{3.11}
\end{equation}

We define $u(t)\triangleq\exp\{\int^t_02k ds\}(c+k^2+\lambda_2^2)$, from
\eqref{3.11} we know
$$
\begin{aligned}
u'(t)&=2k \exp\{\int^t_02k ds\}(c+k^2+\lambda_2^2)+\exp\{\int^t_02k
ds\}(2k k'+2\lambda_2 \lambda_2')\\
&=2k\exp\{\int^t_02k
ds\}(-k'-\lambda_1\lambda_2+2\lambda_2^2)+\exp\{\int^t_02k
ds\}(2k k'+2\lambda_2 \lambda_2')\\
&=2\exp\{\int^t_02k
ds\}(-kk'-k\lambda_1\lambda_2+2k\lambda_2^2+kk'+\lambda_2\lambda_2')=0,
\end{aligned}
$$
which means $u(t)$ is a constant. We have the following lemma:

\begin{lem}\label{lemma3.3}
Let $\lambda_1$ and $\lambda_2$ be two real-valued functions of $t$ with $\lambda_1-2\lambda_2\neq 0$, and let
$k(t)=\frac{\lambda_2'}{\lambda_1-2\lambda_2}$, $u=\exp\{\int^t_02k ds\}(c+k^2+\lambda_2^2)$. If $\lambda_1$ and
$\lambda_2$ satisfy \eqref{3.11} then

{\rm(i)} When $u\neq 0$,
$$g_1(t)=(k-\lambda_2i)\exp\{\int^t_0(k-\lambda_2i)ds\},~g_2(t)=\exp\{\int^t_0[(\lambda_2-\lambda_1)i+k]ds\}$$
are two independent complex-valued solutions of the differential
equation
\begin{equation}
g''+\lambda_1ig'+(\lambda_1'i+c)g=0, \label{3.13}
\end{equation}

{\rm(ii)} When $u=0$,
$\tilde{g}_1(t)=g_2(t)\int^t_0[\exp\{\int^s_0[(\lambda_1-2\lambda_2)i-2k]dx\}]ds$
and $g_2(t)$

are two independent complex-valued solutions of \eqref{3.13}.
\end{lem}
\begin{proof}
We define \begin{equation}
\begin{aligned}
&f_1(t)=-g_1'(t)-\lambda_1ig_1(t),\\
&f_2(t)=-g_2'(t)-\lambda_1ig_2(t),\\
&\tilde{f}_1(t)=-\tilde{g}_1'(t)-\lambda_1i\tilde{g}_1(t).
\label{35}
\end{aligned}
\end{equation}
by using \eqref{3.11}, after a direct calculation we have
\begin{equation}
\begin{aligned}
&f_1(t)=c\exp\{\int^t_0(k-\lambda_2i)ds\},\\
&f_2(t)=-\exp\{\int^t_0[(\lambda_2-\lambda_1)i+k]ds\}(k+\lambda_2i)=-g_2(t)(k+\lambda_2i),\\
&\tilde{f}_1(t)=-(k+\lambda_2i)\tilde{g}_1(t)-g_2(t)\exp\{\int^t_0[(\lambda_1-2\lambda_2)i-2k]ds\},\label{44}
\end{aligned}
\end{equation}
which implies that $f'_1(t)=cg_1(t),~f'_2(t)=cg_2(t)$ and
$\tilde{f}'_1(t)=c\tilde{g}_1(t)$. Hence, we obtain $g_1(t)$,
$g_2(t)$ and $\tilde{g}_1(t)$ are all complex-valued solutions of
\eqref{3.13}.

We denote that
$$
\begin{aligned}
&f(t)=\frac{g_1(t)}{g_2(t)}=(k-\lambda_2 i)\exp\{\int^t_0[(\lambda_1-2\lambda_2)i]ds\},\\
&\tilde{f}(t)=\frac{\tilde{g}_1(t)}{g_2(t)}=\int^t_0[\exp\{\int^s_0[(\lambda_1-2\lambda_2)i-2k]dx\}]ds,
\end{aligned}
$$
by a straightforward calculation, we have
\begin{equation}
f'(t)=-u\exp\{\int^t_0[(\lambda_1-2\lambda_2)i-2k]ds\},~\tilde{f}'(t)=\exp\{\int^t_0[(\lambda_1-2\lambda_2)i-2k]ds\}.\label{48}
\end{equation}

(i) When $u\neq 0$,  by definition we have $g_1(t)=f(t)\cdot
g_2(t)$, since $u\neq0$, from \eqref{48} we have $f'(t)\neq0$, which
implies $g_1(t)$ and $g_2(t)$ are independent.

(ii) When $u=0$, by definition we have
$\tilde{g}_1(t)=\tilde{f}(t)\cdot g_2(t),$  from \eqref{48} we have
$\tilde{f}'(t)\neq 0$ ,which implies $\tilde{g}_1(t)$ and $g_2(t)$
are independent.
\end{proof}

\section{Warped product Lagrangian immersions and Calabi product Lagrangian immersions in $\mathbb{CP}^n$}

Throughout this section, we assume that $\psi:M\to \mathbb{CP}^n(4)$ is a Lagrangian immersion,
$M$ admits two orthogonal
distributions $\mathcal{D}_1$ (of dimension 1, spanned by a unit vector $E_1$) and $\mathcal{D}_2$ (of
dimension $n-1$, spanned by $\{E_2,\ldots, E_n\}$), and there exist local functions $\lambda_1,~\lambda_2$ such that
\eqref{1.1} holds,
and $\nabla$$\lambda_1$ lies in $\mathcal{D}_1$.

We use the same notations as in the previous sections.
 We consider a
horizontal lift $\tilde{\psi}$ of $\psi$, and identify
$E_1$ with $\tilde{\psi}_{*}E_1$. We define $\tilde{\phi}_i:M\to
\mathbb{C}^{n+1},i=1,2$, by
$$
\tilde{\phi}_1=f_1\tilde{\psi}+g_1E_1,~
\tilde{\phi}_2=f_2\tilde{\psi}+g_2E_1,
$$
where $f_1$, $g_1$, $f_2$ and $g_2$ are defined in Lemma
\ref{lemma3.3} and \eqref{35}.

From the proof of Lemma \ref{lemma3.3}, it is not difficult to see
that
\begin{equation}
\left\{
\begin{aligned}
&E_1(\tilde{\phi}_1)=E_1(\tilde{\phi}_2)=E_i(\tilde{\phi}_2)=0,~2\leq i\leq n,\\
& \langle \tilde{\phi}_1,\tilde{\phi}_2 \rangle = \langle \tilde{\phi}_1,J\tilde{\phi}_2 \rangle =0,
\label{4.3}
\end{aligned}
\right.
\end{equation}
which   implies that $\tilde{\phi}_2$ is a constant vector. We have
$$
 \langle \tilde{\phi}_1,\tilde{\phi}_1 \rangle = \langle \tilde{\phi}_2,\tilde{\phi}_2 \rangle =\exp\{\int^t_0(2k)ds\}(1+k^2+\lambda_2^2)=u>0.
$$

In particular, let
$\tilde{\psi}_i=\frac{\tilde{\phi}_i}{\sqrt{u}},~i=1,2$, since
$\tilde{\psi}_2$ is a constant unit vector, up to a rigid motion of
$\mathbb{S}^{2n+1}(1)$, we recover $\tilde{\psi}$  in terms of
$\tilde{\psi}_i$, $i=1,2$, by
$$
\begin{aligned}
\tilde{\psi}(t,p)&=-\frac{g_2}{f_2g_1-f_1g_2}\tilde{\psi}_1+\frac{g_1}{f_2g_1-f_1g_2}\tilde{\psi}_2\\
&=(\frac{\exp\{\int^t_0(\lambda_2i+k)ds\}}{\sqrt{u}}\tilde{\psi}_1(0,p),
\frac{(\lambda_2i-k)\exp\{\int^t_0(k+\lambda_1i-\lambda_2i)ds\}}{\sqrt{u}})
\end{aligned}$$
for any $(t,p)\in I\times M_1$.

The properties of the lift $\tilde{\psi}$ imply that
$\tilde{\psi}_1$ is a horizontal immersion in the corresponding
$\mathbb{S}^{2n-1}(1)\subset \mathbb{C}^{n-1}$, which means that
$\tilde{\psi}_1$ is the horizontal lift of a Lagrangian immersion
$\psi_1=\Pi(\tilde{\psi}_1),$ in $\mathbb{CP}^{n-1}(4)$. It is not
difficult to check that
$\tilde{\gamma}(t)=(\frac{\exp\{\int^t_0(\lambda_2i+k)ds\}}{\sqrt{u}},\frac{(\lambda_2i-k)\exp\{\int^t_0(k+\lambda_1i-\lambda_2i)ds\}}{\sqrt{u}})$
satisfies the horizontality condition and hence is a Legendre curve
in $\mathbb{S}^3(1)$, so we have $\psi$ is the warped product
Lagrangian immersion of a point with a (n-1)-dimensional Lagrangian
submanifold $\psi_1$.

Hence we have proved the following  theorem:
\begin{thm}\label{thm4.1}
Let $\psi:M\to \mathbb{CP}^n(4)(n\geq3)$ be a Lagrangian immersion.
Assume $M$ admits two orthogonal distributions $\mathcal{D}_1$ (of
dimension 1, spanned by a unit vector $E_1$) and $\mathcal{D}_2$ (of
dimension $n-1$, spanned by $\{E_2,\ldots, E_n\}$), and there exist
local functions $\lambda_1,~\lambda_2$ such that \eqref{1.1} holds.
Moreover we assume that $\nabla$$\lambda_1$ lies in $\mathcal{D}_1$
and let
$k=\frac{E_1(\lambda_2)}{\lambda_1-2\lambda_2},~u=\exp\{\int^t_02k
ds\}(1+k^2+\lambda_2^2)$. We have that $M$ is locally isometric to a
product $I\times M_1$ and $\psi$ is locally a warped product
Lagrangian immersion of a point with a (n-1)-dimensional Lagrangian
submanifold $\psi_1: M_1\to \mathbb{CP}^{n-1}(4)$, and up to a rigid
motion of $\mathbb{CP}^{n}(4)$, $\psi: I\times M_1\to
\mathbb{CP}^n(4)$ is locally given by $\psi=\Pi\circ\tilde{\psi}$
with
$$
\tilde{\psi}(t,p)=(\frac{\exp\{\int^t_0(\lambda_2i+k)ds\}}{\sqrt{u}}\tilde{\psi}_1(p),
\frac{(\lambda_2i-k)\exp\{\int^t_0(k+\lambda_1i-\lambda_2i)ds\}}{\sqrt{u}}),
$$
where $\Pi$ is the Hopf fibration and $\tilde{\psi}_1:M_1\to
\mathbb{S}^{2n-1}(1)$ is the horizontal lift of $\psi_1$.
\end{thm}

\begin{rem}\label{rem42}
In Theorem \ref{thm4.1}, if we assume that both $\nabla\lambda_1$
and $\nabla\lambda_2$ lie in $\mathcal{D}_1$, we can get the same
conclusion without assuming $n\geq3$.
\end{rem}

By applying Theorem \ref{thm4.1}, we can  prove Theorem
\ref{thm1.3}.

\medskip
\noindent{\bf Proof of Theorem \ref{thm1.3} :} Suppose $M$ admits
two orthogonal distributions $\mathcal{D}_1$ (of dimension 1,
spanned by a unit vector $E_1$) and $\mathcal{D}_2$ (of dimension
$n-1$, spanned by $\{E_2,\ldots, E_n\}$), and there exist two real
constants $\lambda_1$ and $\lambda_2$ such that \eqref{1.1} holds.
Since $\lambda_1$ and $\lambda_2$ are both constants, we have $k=0$,
$u=1+\lambda_2^2$ and $\lambda_1\lambda_2-\lambda_2^2+1=0$. Hence
from Theorem \ref{thm4.1} and Remark \ref{rem42}, we have $M$ is
locally isometric to a product $I\times M_1$ and  up to a
reparametrization and a rigid motion of $\mathbb{CP}^{n}(4)$, $\psi:
I\times M_1\to \mathbb{CP}^n(4)$ is locally given by
$\psi=\Pi\circ\tilde{\psi}$ with
\begin{equation}
\tilde{\psi}(t,p)=(\frac{1}{\sqrt{\lambda_2^2+1}}e^{\frac{it\lambda_2^2}{1+\lambda_2^2}}\tilde{\psi}_1(p),
\frac{\lambda_2}{\sqrt{\lambda_2^2+1}}e^{\frac{-it}{1+\lambda_2^2}}),~(t,p)\in
I\times M_1.
\end{equation}
By Definition \ref{def1.2} we conclude that $\psi$ is
locally a Calabi product Lagrangian immersion of an (n-1)-dimensional Lagrangian
immersion $\psi_1: M_1\to \mathbb{CP}^{n-1}(4)$ and a point.
By applying Corollary 4.1 in \cite{CLU}, we have

{\rm(i)}\quad $\psi$ is Hamiltonian minimal if and only if $\psi_1$
is Hamiltonian minimal.

{\rm(i)}\quad $\psi$ is minimal if and only if $\lambda_2=\pm
\frac{1}{\sqrt{n}}$ and $\psi_1$ is minimal. In this case, up to a reparametrization and a
rigid motion of $\mathbb{CP}^n$, $\psi: I\times M_1\to
\mathbb{CP}^n(4)$ is locally given by $\psi=\Pi\circ\tilde{\psi}$ with
\begin{equation}
\tilde{\psi}(t,p)=(\sqrt{\frac{n}{n+1}}e^{i\frac{1}{n+1}t}\tilde{\psi}_1(p),
\sqrt{\frac{1}{n+1}}e^{-i\frac{n}{n+1}t}),~(t,p)\in I\times
M_1,\label{4-9}
\end{equation}
where $\Pi$ is the Hopf fibration and $\tilde{\psi}_1:M_1\to
\mathbb{S}^{2n-1}(1)$ is the horizontal lift of $\psi_1$.

Conversely,  if $\psi$ is locally a Calabi product Lagrangian
immersion of an (n-1)-dimensional Lagrangian immersion $\psi_1:
M_1\to \mathbb{CP}^{n-1}(4)$ and a point, by Definition \ref{def1.2}
we assume that
$\psi(t,p)=\Pi(r_1e^{i(\frac{r_2}{r_1}at)}\tilde{\psi}_1(p),r_2e^{i(-\frac{r_1}{r_2}at)}),$
where $r_1$, $r_2$ and $a$ are positive constants,
$\tilde{\psi}_1:M_1\to \mathbb{S}^{2n-1}(1)$ is the horizontal lift
of $\psi_1$. For tangent vectors $e_i\in M_1$  we define local
vector fields
$$E_1=\frac{\psi_t}{|\psi_t|},~E_i=\frac{\psi_*(0,e_i)}{|\psi_*(0,e_i)|},$$
after a straightforward calculation we have the second fundamental
form satisfies
$$h(E_1,E_1)=(\frac{r_2}{r_1}-\frac{r_1}{r_2})JE_1,~h(E_1,E_i)=\frac{r_2}{r_1}JE_i,~(\frac{r_2}{r_1}-\frac{r_1}{r_2})\neq 2\frac{r_2}{r_1}.$$

Hence, we have completed the proof of of Theorem \ref{thm1.3}.\qed

\begin{rem}\label{rem4.2}
In \eqref{4-9}, if $\psi_1$ is a totally geodesic Lagrangian
embedding of $\mathbb{S}^{n-1}$ into $\mathbb{CP}^{n-1}(4)$, then we
get a minimal Lagrangian embedding
$$\overline{(e^{it},x)}\mapsto\Pi(\sqrt{\frac{n}{n+1}}e^{i\frac{1}{n+1}t}x,
\sqrt{\frac{1}{n+1}}e^{-i\frac{n}{n+1}t})$$
of the quotient $(\mathbb{S}^1\times\mathbb{S}^{n-1})/\mathbb{Z}_2$ into $\mathbb{CP}^n$, where the action of $\mathbb{Z}$ is generated
by the involutions $(e^{it},x)\mapsto (-e^{it},-x),~(e^{it},x)\mapsto (-e^{it},x).$ This example has been studied in \cite{N} and \cite{CLU}.
\end{rem}

\noindent{\bf Proof of Theorem \ref{thm1.4}:}
Suppose $M$ has parallel second fundamental form,  $M$ admits two orthogonal
distributions $\mathcal{D}_1$ (of dimension 1, spanned by a unit vector $E_1$) and $\mathcal{D}_2$ (of
dimension $n-1$, spanned by $\{E_2,\ldots, E_n\}$), and there exist local functions $\lambda_1,~\lambda_2$ such that
\eqref{1.1} holds.

By
\eqref{23} and \eqref{1.1} we have
\begin{equation}
\left\{
\begin{aligned}
(\nabla h)(E_1,E_1,E_1)&=E_1(\lambda_1)JE_1+(\lambda_1-2\lambda_2)J\nabla_{E_1}E_1,\\
(\nabla h)(E_i,E_1,E_1)&=E_i(\lambda_1)JE_1+(\lambda_1-2\lambda_2)\nabla_{E_i}E_1,\\
(\nabla h)(E_1,E_i,E_1)&=E_1(\lambda_2)JE_i-h(\nabla_{E_1}E_1,E_i)\\
&+(\lambda_2-\lambda_1) \langle \nabla_{E_1}E_i,E_1 \rangle JE_1,\\
(\nabla
h)(E_i,E_j,E_1)&=E_i(\lambda_2)JE_j-h(\nabla_{E_i}E_1,E_j)\\
&+(\lambda_1-\lambda_2) \langle \nabla_{E_i}E_1,E_j \rangle JE_1.
\label{4.6}
\end{aligned}
\right.
\end{equation}
Since $M$ has parallel second fundamental form and
$\lambda_1\neq2\lambda_2$, \eqref{4.6} implies that
$$
E_1(\lambda_1)=E_i(\lambda_1)=E_1(\lambda_2)=E_i(\lambda_2)=
\nabla_{E_1}E_1=\nabla_{E_i}E_1=0. $$
Hence we have $\lambda_1$ and $\lambda_2$ are both real constants,  by applying Theorem \ref{thm1.3} we have $\psi$ is
locally a Calabi product Lagrangian immersion of a point and an (n-1)-dimensional Lagrangian
immersion $\psi_1: M_1\to \mathbb{CP}^{n-1}(4)$.

A straightforward calculation gives the following claim:
\begin{cla}\label{cla4.2}
a Calabi product Lagrangian immersion of
an (n-1)-dimensional Lagrangian
immersion $\psi_1: M_1\to \mathbb{CP}^{n-1}(4)$ and a point has parallel second fundamental form if and only if $\psi_1$ has parallel
second fundamental form.
\end{cla}

In fact, suppose $\psi$ is
a Calabi product Lagrangian immersion of an (n-1)-dimensional Lagrangian
immersion $\psi_1: M_1\to \mathbb{CP}^{n-1}(4)$ and a point, by Definition \ref{def1.2} we assume that
$$\psi(t,p)=\Pi(r_1e^{i(\frac{r_2}{r_1}at)}\tilde{\psi}_1(p),r_2e^{i(-\frac{r_1}{r_2}at)}),$$
where $r_1$, $r_2$ and $a$ are positive constants, $\tilde{\psi}_1:M_1\to
\mathbb{S}^{2n-1}(1)$ is the horizontal lift of $\psi_1$. Let $\nabla^1$ denote the induced connection of $M_1$ and $h^1$ denote the second
fundamental form of $\psi_1: M_1\to \mathbb{CP}^{n-1}(4)$, we choose $E_1$ and $E_i$ as in the proof of Theorem \ref{thm1.3},
after a straightforward computation, we have
\begin{equation}
\begin{aligned}
&\nabla_{E_1}E_1=\nabla_{E_i}E_1=0,~\nabla_{E_i}E_j=\nabla^1_{E_i}E_j,\\
&h(E_1,E_1)=(\frac{r_2}{r_1}-\frac{r_1}{r_2})JE_1,~h(E_1,E_i)=\frac{r_2}{r_1}JE_i,\\
&h(E_i,E_j)=h^1(E_i,E_j)+\frac{r_2}{r_1}\delta_{ij}JE_1.
\end{aligned}\label{4.12}
\end{equation}
From \eqref{23} and \eqref{4.12} we have
$$
(\nabla h)(E_i, E_j,E_k)=(\nabla ^1h^1)(E_i, E_j,E_k)
$$
and all the other components of $\nabla h$ are zero. Therefore, we have proved Claim \ref{cla4.2} and
completed the proof of Theorem \ref{thm1.4}.\qed

After an analogous argument as in the proof of Theorem \ref{thm1.3} and Theorem \ref{thm1.4}, by using
Theorem \ref{thm4.4} of \cite{BMV}, we have

\begin{thm}\label{thm4.3}
Let $\psi:M\to \mathbb{CP}^n(4)$ be a Lagrangian immersion, $\psi$ is locally a Calabi product Lagrangian immersion of two lower dimensional Lagrangian
immersion $\psi_i: M_1\to \mathbb{CP}^{n_1}(4)$ and $\psi_2: M_2\to \mathbb{CP}^{n_2}(4)$
if and only if $M$ admits three mutually orthogonal
distributions $\mathcal{D}_1$ (spanned by a unit vector $E_1$),
$\mathcal{D}_2$, and $\mathcal{D}_3$ of dimension $1$, $n_1$ and $n_2$ respectively, with
$1+n_1+n_2=n$, and there three real constants
$\lambda_1$, $\lambda_2$ and $\lambda_3$ ($2\lambda_3\neq\lambda_1\neq2\lambda_2\neq2\lambda_3$) s.t.
for all $E_i\in\mathcal{D}_2,~E_\alpha\in\mathcal{D}_3$,
\begin{equation}
\left\{
\begin{aligned}
&h(E_1,E_1)=\lambda_1JE_1,~h(E_1,E_i)=\lambda_2JE_i,\\
&h(E_1,E_\alpha)=\lambda_3JE_\alpha,~h(E_i,E_\alpha)=0.
\end{aligned}\label{4.16}
\right.\end{equation}
Moreover, a Lagrangian immersion  $\psi:M\to \mathbb{CP}^n(4)$ satisfying the above conditions has the following properties:

{\rm(i)}\quad $\psi$ is Hamiltonian minimal if  both $\psi_1$ and
$\psi_2$ are Hamiltonian minimal.

{\rm(ii)}\quad $\psi$ is minimal if and only if both $\psi_1$ and $\psi_2$ are minimal and
$$\lambda_1=\pm\frac{n_2-n_1}{\sqrt{(n_1+1)(n_2+1)}},~\lambda_2=\pm\frac{\sqrt{n_2+1}}{\sqrt{n_1+1}},~\lambda_3=\mp\frac{\sqrt{n_1+1}}{\sqrt{n_2+1}}.$$
In this case, up to a reparametrization and a
rigid motion of $\mathbb{CP}^n$,
$\psi$ is locally given by
$\psi=\Pi\circ\tilde{\psi}$ with
$$
\tilde{\psi}=(\sqrt{\frac{n_1+1}{n_1+n_2+2}}e^{i\frac{n_2+1}{n_1+n_2+2}t}\tilde{\psi}_1,
\sqrt{\frac{n_2+1}{n_1+n_2+2}}e^{-i\frac{n_1+1}{n_1+n_2+2}t}\tilde{\psi}_2),
$$
where $\Pi$ is the Hopf fibration,  $\tilde{\psi}_1:M_1\to
\mathbb{S}^{2n_1-1}(1)$ and $\tilde{\psi}_2:M_2\to
\mathbb{S}^{2n_2-1}(1)$ are horizontal lifts of $\psi_1$ and
$\psi_2$.
\end{thm}

\begin{thm}\label{thm4.4}
Let $\psi:M\to \mathbb{CP}^n(4)$ be a Lagrangian immersion. If $M$ admits three mutually orthogonal
distributions $\mathcal{D}_1$ (spanned by a unit vector $E_1$),
$\mathcal{D}_2$, and $\mathcal{D}_3$ of dimension $1$, $n_1$ and $n_2$ respectively, with
$1+n_1+n_2=n$, and there three real constants
$\lambda_1$, $\lambda_2$ and $\lambda_3$ ($2\lambda_3\neq\lambda_1\neq2\lambda_2\neq2\lambda_3$) such that
\eqref{4.16} holds for all $E_i\in\mathcal{D}_2,~E_\alpha\in\mathcal{D}_3$,
then $M$ has parallel second fundamental form if and only if $\psi$ is locally a Calabi product Lagrangian immersion of two lower dimensional Lagrangian submanifolds
$\psi_i(i=1,2)$ with parallel second fundamental form.
\end{thm}

\section{Warped product Lagrangian immersions and Calabi product Lagrangian immersions in $\mathbb{CH}^n$}
Throughout this section, we assume that $\psi:M\to
\mathbb{CH}^n(-4)$ is a Lagrangian immersion, $M$ admits two
orthogonal distributions $\mathcal{D}_1$ (of dimension 1, spanned by
a unit vector $E_1$) and $\mathcal{D}_2$ (of dimension $n-1$,
spanned by $\{E_2,\ldots, E_n\}$), and there exist local functions
$\lambda_1,~\lambda_2$ such that \eqref{1.1} holds, and
$\nabla$$\lambda_1$ lies in $\mathcal{D}_1$.

We use the same notations as in the previous sections.
 We consider a
horizontal lift $\tilde{\psi}$ of $\psi$, and identify $E_1$ with
$\tilde{\psi}_{*}E_1$.
 We define $\tilde{\phi}_i:M\to
\mathbb{C}_1^{n+1},i=0,1,2$, by
$$
\tilde{\phi}_0=\tilde{f}_1\tilde{\psi}+\tilde{g}_1E_1,~
\tilde{\phi}_1=f_1\tilde{\psi}+g_1E_1,~
\tilde{\phi}_2=f_2\tilde{\psi}+g_2E_1,
$$
where $\tilde{f}_1$,  $\tilde{g}_1$, $f_1$, $g_1$, $f_2$ and $g_2$
and are defined in Lemma \ref{lemma3.3} and \eqref{35}.

In section 3, we have proved that
$u=\exp\{\int^t_0(2k)ds\}(-1+k^2+\lambda_2^2)$ is a constant, when
$u\neq 0$, from the proof of Lemma \ref{lemma3.3}, it is not
difficult to see that
\begin{equation}
\left\{
\begin{aligned}
&E_1(\tilde{\phi}_1)=E_1(\tilde{\phi}_2)=E_i(\tilde{\phi}_2)=0,~2\leq i\leq n,\\
& \langle \tilde{\phi}_1,\tilde{\phi}_2 \rangle_1 = \langle \tilde{\phi}_1,J\tilde{\phi}_2 \rangle_1 =0.
\label{5.3}
\end{aligned}
\right.
\end{equation}
By \eqref{5.3} we  have $\tilde{\phi}_2$ is a constant vector. We
have
$$
 \langle \tilde{\phi}_1,\tilde{\phi}_1 \rangle_1 =- \langle \tilde{\phi}_2,\tilde{\phi}_2 \rangle_1 =\exp\{\int^t_0(2k)ds\}(-1+k^2+\lambda_2^2)=u.
$$

We have three cases:

{\rm(i)}\quad $u>0$, let
$\tilde{\psi}_i=\frac{\tilde{\phi}_i}{\sqrt{u}},~i=1,2$, since
$\tilde{\psi}_2$ is a constant vector with $
\langle\tilde{\psi}_2,\tilde{\psi}_2\rangle_1=-1$, up to a rigid
motion of $\mathbb{H}_1^{2n+1}(-1)$, we recover $\tilde{\psi}$ in
terms of $\tilde{\psi}_i$, $i=1,2$,  by
$$
\begin{aligned}
\tilde{\psi}(t,p)&=-\frac{g_2}{f_2g_1-f_1g_2}\tilde{\psi}_1+\frac{g_1}{f_2g_1-f_1g_2}\tilde{\psi}_2\\
&=(\frac{(\lambda_2i-k)\exp\{\int^t_0(k+\lambda_1i-\lambda_2i)ds\}}{\sqrt{u}},\frac{\exp\{\int^t_0(\lambda_2i+k)ds\}}{\sqrt{u}}\tilde{\psi}_1(0,p))
\end{aligned}
$$
for any $(t,p)\in I\times M_1$.

{\rm(ii)}\quad $u<0$, let
$\tilde{\psi}_i=\frac{\tilde{\phi}_i}{\sqrt{-u}},~i=1,2$, since
$\tilde{\psi}_2$ is a constant vector with $
\langle\tilde{\psi}_2,\tilde{\psi}_2\rangle_1=1$,  up to a rigid
motion of $\mathbb{H}_1^{2n+1}(-1)$, we recover $\tilde{\phi}$ in
terms of $\tilde{\psi}_i$, $i=1,2$,  by
$$
\begin{aligned}
\tilde{\psi}(t,p)&=-\frac{g_2}{f_2g_1-f_1g_2}\tilde{\psi}_1+\frac{g_1}{f_2g_1-f_1g_2}\tilde{\psi}_2\\
&=(\frac{\exp\{\int^t_0(\lambda_2i+k)ds\}}{\sqrt{-u}}\tilde{\psi}_1(0,p),\frac{(\lambda_2i-k)\exp\{\int^t_0(k+\lambda_1i-\lambda_2i)ds\}}{\sqrt{-u}})
\end{aligned}
$$
for any $(t,p)\in I\times M_1$.

For case (i) and case (ii), an analogous reasoning to that in
Theorem \ref{thm4.1} yields that $\psi$ is the warped product
Lagrangian immersion of a point with a (n-1)-dimensional Lagrangian
submanifold $\psi_1$.

{\rm(iii)}\quad $u=0$, by \eqref{48} we have $ f(t)=(k-\lambda_2
i)\exp\{\int^t_0[(\lambda_1-2\lambda_2)i]ds\} $ is a constant
function. Hence we have $ f(t)\equiv k(0)-\lambda_2(0)i. $ By a
direct calculation, we have $
f_2\tilde{g_1}-\tilde{f_1}g_2=\exp\{\int_0^t(-\lambda_1i)ds\}$ and
$$
\begin{aligned}
&\tilde{\psi}=e^{\int_o^t(k+\lambda_2i)ds}\Big((k(0)-\lambda_2(0)i)\cdot
\int_0^t[(k+\lambda_2i)e^{\int_0^s-2kdx}]ds\cdot \tilde{\psi}_2-\tilde{\psi}_0\Big),\\
& \langle \tilde{\psi}_2,\tilde{\psi}_0 \rangle_1 =-k(0),~ \langle
\tilde{\psi}_2,J\tilde{\psi}_0 \rangle_1 =-\lambda_2(0).
\end{aligned}$$

We choose the initial conditions at the origin $0\in
\mathbb{C}_1^{n+1}$ such that
$$
\tilde{\psi}(0)=(1,0,\ldots,0),~\tilde{\psi}_t(0)=(0,1,0,\ldots,0),
$$
and denote that $\tilde{\psi}_0=(-A_0,-A_1,-\psi_3)$, where $A_0$
and $A_1$ are two complex-valued functions. By using the fact that $
\langle \tilde{\psi},\tilde{\psi} \rangle _1=-1$, we obtain
\begin{equation}
\begin{aligned}&\tilde{\psi}_2=(-k(0)-\lambda_2(0)i,1,0,\ldots,0),\\
&A_1=(i\lambda_2(0)-k(0))(A_0-1),Re(A_0)=1+\frac{ \langle \psi_3,\psi_3 \rangle }{2},
\end{aligned}\end{equation}
and
\begin{equation}
\begin{aligned}
\tilde{\psi}(t,u_2,\ldots,u_n)&=e^{\int_o^t(k+\lambda_2i)ds}\bigg(-\int_0^t(k+\lambda_2i)e^{\int_0^s-2kdx}ds+A_0,\\
&(i\lambda_2(0)-k(0))
(-\int_0^t(k+\lambda_2i)e^{\int_0^s-2kdx}ds+A_0-1),\psi_3\bigg).
\end{aligned}
\label{5.12}
\end{equation}

The immersion  $\tilde{\psi}$ given by \eqref{5.12} is a Lagrangian
immersion in $\mathbb{CH}^n(-4)$ if and only if $\psi_3:M_1\to
\mathbb{C}^{n-1}$ is Lagrangian and $A_0$ satisfies the following
conditions:
$$
Re(A_0)=1+\frac{ \langle \psi_3,\psi_3 \rangle }{2},~v(Im
A_0)= \langle \psi_{3*}(v),J\psi_3 \rangle ,~\forall v\in T_pM_1, $$

Hence we have proved the following theorem:

\begin{thm}\label{thm5.1}
Let $\psi:M\to \mathbb{CH}^n(-4)(n\geq3)$ be a Lagrangian immersion.
Assume $M$ admits two orthogonal
distributions $\mathcal{D}_1$(of dimension 1) and $\mathcal{D}_2$(of
dimension $n-1$) such that \eqref{1.1} holds. Moreover we assume that $\nabla$$\lambda_1$ lies in $\mathcal{D}_1$ and let
$k=\frac{E_1(\lambda_2)}{\lambda_1-2\lambda_2},~u=\exp\{\int^t_02k
ds\}(-1+k^2+\lambda_2^2)$. We have

(1) If  $u>0$, then $M$ is locally isometric to a product $I\times
M_1$ and $\psi$ is locally a warped product Lagrangian immersion of
a point with a (n-1)-dimensional Lagrangian submanifold $\psi_1:
M_1\to \mathbb{CP}^{n-1}(4)$, and up to a rigid motion of
$\mathbb{CH}^{n}(-4)$, $\psi: I\times M_1\to \mathbb{CH}^n(-4)$ is
locally given by $\psi=\Pi\circ\tilde{\psi}$ with
$$
\tilde{\psi}(t,p)=(\frac{(\lambda_2i-k)\exp\{\int^t_0(k+\lambda_1i-\lambda_2i)ds\}}{\sqrt{u}},\frac{\exp\{\int^t_0(\lambda_2i+k)ds\}}
{\sqrt{u}}\tilde{\psi}_1(0,p)), $$
where $\Pi$ is the Hopf fibration and $\tilde{\psi}_1:M_1\to
\mathbb{S}^{2n-1}(1)$ is the horizontal lift of $\psi_1$.

(2) If $u<0$, then $M$ is locally isometric to a product $I\times
M_1$ and $\psi$ is locally a warped product Lagrangian immersion of
a point with a (n-1)-dimensional Lagrangian submanifold $\psi_1:
M_1\to \mathbb{CH}^{n-1}(-4)$, and and up to a rigid motion of
$\mathbb{CH}^{n}(-4)$, $\psi: I\times M_1\to \mathbb{CH}^n(-4)$ is
locally given by $\psi=\Pi\circ\tilde{\psi}$ with
$$
\tilde{\psi}(t,p)=(\frac{\exp\{\int^t_0(\lambda_2i+k)ds\}}{\sqrt{-u}}\tilde{\psi}_1(0,p),\frac{(\lambda_2i-k)
\exp\{\int^t_0(k+\lambda_1i-\lambda_2i)ds\}}{\sqrt{-u}}),$$
where $\Pi$ is the Hopf fibration and $\tilde{\psi}_1:M_1\to
\mathbb{H}_1^{2n-1}(-1)$ is the horizontal lift of $\psi_1$.

(3) If $u=0$, then $M$ is locally isometric to a product $I\times
M_1$ and up to a rigid motion of $\mathbb{CH}^n(-4)$, $\psi$ is
locally given by $\Pi\circ\tilde{\psi}$, where $\Pi$ is the Hopf
fibration and $\tilde{\psi}:I\times M^{n-1}\to
\mathbb{H}^{2n+1}_1(-1)$ is given by
$$
\begin{aligned}
\tilde{\psi}(t,u_2,\ldots,u_n)&=e^{\int_o^t(k+\lambda_2i)ds}\bigg(-\int_0^t(k+\lambda_2i)e^{\int_0^s-2kdx}ds+A_0,\\
&(i\lambda_2(0)-k(0))
(-\int_0^t(k+\lambda_2i)e^{\int_0^s-2kdx}ds+A_0-1),\psi_3\bigg),
\end{aligned}
$$
where $A_0$ is a complex function on $M_1$ and $\psi_3:M_1\to
\mathbb{C}^{n-1}$ is a Lagrangian immersion, and they
satisfy that $$Re A_0=\frac{1}{2} \langle \psi_3,\psi_3 \rangle +1;~v(Im
A_0)= \langle \psi_{3*}(v),J\psi_3 \rangle ,~\forall v\in T_pM_1.$$
\end{thm}

\begin{rem}\label{rem5.2}
In Theorem \ref{thm5.1}, if we assume that both $\nabla\lambda_1$
and $\nabla\lambda_2$ lie in $\mathcal{D}_1$, we can get the same
conclusion without assuming $n\geq3$.
\end{rem}

By applying Theorem \ref{thm5.1}, we can  prove Theorem
\ref{thm1.5}.

\medskip
 \noindent{\bf Proof of Theorem  \ref{thm1.5} :} Suppose
$M$ admits two orthogonal distributions $\mathcal{D}_1$ (of
dimension 1, spanned by a unit vector $E_1$) and $\mathcal{D}_2$ (of
dimension $n-1$, spanned by $\{E_2,\ldots, E_n\}$), and there exist
two real constants $\lambda_1$ and $\lambda_2$  such that
\eqref{1.1} holds. Since $\lambda_1$ and $\lambda_2$ are both
constants, we have $k=0$, $u=-1+\lambda_2^2$ and
$\lambda_1\lambda_2-\lambda_2^2-1=0$, $\lambda_1\neq 2\lambda_2$.
 Hence $\lambda_1=\lambda_2+\frac{1}{\lambda_2}\neq 2\lambda_2$, which implies that $u=-1+\lambda_2^2\neq 0$ and
 $\lambda_1+(n-1)\lambda_2=\frac{n\lambda_2^2+1}{\lambda_2}\neq 0$, so $M$ can never be minimal.
From Theorem \ref{thm5.1} and Remark \ref{rem5.2}, we have

(1) If $\lambda_2^2-1>0$, then $M$ is locally isometric to a product
$I\times M_1$ and  up to a reparametrization and a rigid motion of
$\mathbb{CH}^{n}(-4)$,  $\psi: I\times M_1\to \mathbb{CH}^n(-4)$ is
locally given by $\psi=\Pi\circ\tilde{\psi}$ with
$$
\tilde{\psi}(t,p)=(\frac{\lambda_2}{\sqrt{\lambda_2^2-1}}e^{\frac{it}{\lambda_2^2-1}},
\frac{1}{\sqrt{\lambda_2^2-1}}e^{\frac{it\lambda_2^2}{\lambda_2^2-1}}\tilde{\psi}_1(p)),~(t,p)\in
I\times M_1,
$$
where $\Pi$ is the Hopf fibration and $\tilde{\psi}_1:M_1\to
\mathbb{S}^{2n-1}(1)$ is the horizontal lift of $\psi_1$.

(2) If $\lambda_2^2-1<0$, then $M$ is locally isometric to a product
$I\times M_1$ and up to a reparametrization and a rigid motion of
$\mathbb{CH}^{n}(-4)$,  $\psi: I\times M_1\to \mathbb{CH}^n(-4)$ is
locally given by $\psi=\Pi\circ\tilde{\psi}$ with
$$
\tilde{\psi}(t,p)=(\frac{1}{\sqrt{1-\lambda_2^2}}e^{\frac{it\lambda_2^2}{1-\lambda_2^2}}\tilde{\psi}_1(p),
\frac{\lambda_2}{\sqrt{1-\lambda_2^2}}e^{\frac{it}{1-\lambda_2^2}}),~(t,p)\in
I\times M_1,
$$
where $\Pi$ is the Hopf fibration and $\tilde{\psi}_1:M_1\to
\mathbb{H}_1^{2n-1}(-1)$ is the horizontal lift of $\psi_1$.

For both Case (1) and Case (2), by Definition \ref{def1.3} we
conclude that $\psi$ is locally a Calabi product Lagrangian
immersion of an (n-1)-dimensional Lagrangian immersion $\psi_1:
M_1\to \mathbb{CP}^{n-1}(4)$ (or $\psi_1: M_1\to
\mathbb{CH}^{n-1}(-4)$) and a point. By applying Corollary 6.6 in
\cite{CLU}, we have $\psi$ is Hamiltonian minimal if and only if
$\psi_1$ is Hamiltonian minimal.

Conversely,  if $\psi$ is
locally a Calabi product Lagrangian immersion of an (n-1)-dimensional Lagrangian
immersion and a point, by Definition \ref{def1.3} we
have two cases:

(1) $
\psi=\Pi(r_1e^{i(\frac{r_2}{r_1}at)},r_2e^{i(-\frac{r_1}{r_2}at)}\tilde{\psi}_1),$
 where $r_1$, $r_2$ and $a$ are positive constants, $\tilde{\psi}_1:M_1\to
\mathbb{S}^{2n-1}(1)$ is the horizontal lift of a (n-1)-dimensional
Lagrangian immersion $\psi_1:M_1\to \mathbb{CP}^{n-1}(4)$. For
tangent vectors $e_i\in M_1$  we can define local vector fields
$$E_1=\frac{\psi_t}{|\psi_t|},~E_i=\frac{\psi_*(0,e_i)}{|\psi_*(0,e_i)|}.$$
It is not difficult to verify that the second fundamental form
satisfies
$$h(E_1,E_1)=(\frac{r_2}{r_1}+\frac{r_1}{r_2})JE_1,~h(E_1,E_i)=\frac{r_1}{r_2}JE_i.$$

 (2) $
\psi=\Pi(r_1e^{i(\frac{r_2}{r_1}at)}\tilde{\psi}_1,r_2e^{i(-\frac{r_1}{r_2}at)}),$
 where $r_1$, $r_2$ and $a$ are positive constants, $\tilde{\psi}_1:M_1\to
\mathbb{H}_1^{2n-1}(1)$ is the horizontal lift of a
(n-1)-dimensional Lagrangian immersion $\psi_1:M_1\to
\mathbb{CH}^{n-1}(-4)$. For tangent vectors $e_i\in M_1$  we define
local vector fields
$$E_1=\frac{\psi_t}{|\psi_t|},~E_i=\frac{\psi_*(0,e_i)}{|\psi_*(0,e_i)|}.$$
It is not difficult to verify that the second fundamental form satisfies
$$h(E_1,E_1)=(\frac{r_2}{r_1}+\frac{r_1}{r_2})JE_1,~h(E_1,E_i)=\frac{r_2}{r_1}JE_i.$$

In both case (1) and case (2), since $-r_1^2+r_2^2=-1$, we have
$(\frac{r_2}{r_1}+\frac{r_1}{r_2})\neq 2 \frac{r_1}{r_2}$ and
$(\frac{r_2}{r_1}+\frac{r_1}{r_2})\neq 2 \frac{r_2}{r_1}.$ Hence, we
have completed the proof of of Theorem \ref{thm1.5}.\qed

After an analogous argument as in the proof of Theorem 4.4 in \cite{BMV} and in the proof of Theorem \ref{thm1.5}, we have the following
result.
\begin{thm}\label{thm5.2}
Let $\psi:M\to \mathbb{CH}^n(-4)$ be a Lagrangian immersion. $\psi$ is locally a Calabi product Lagrangian immersion of two lower dimensional Lagrangian
immersion $\psi_1: M_1\to \mathbb{CP}^{n_1}(4)$ and $\psi_2: M_2\to \mathbb{CH}^{n_2}(-4)$ if and only if
$M$ admits three mutually orthogonal
distributions $\mathcal{D}_1$ (spanned by a unit vector $E_1$),
$\mathcal{D}_2$, and $\mathcal{D}_3$ of dimension $1$, $n_1$ and $n_2$ respectively, with
$1+n_1+n_2=n$, and there three real constants
$\lambda_1$, $\lambda_2$ and $\lambda_3$ ($2\lambda_3\neq\lambda_1\neq2\lambda_2\neq2\lambda_3$) s.t.
for all $E_i\in\mathcal{D}_2,~E_\alpha\in\mathcal{D}_3$,
$$
\left\{
\begin{aligned}
&h(E_1,E_1)=\lambda_1JE_1,~h(E_1,E_i)=\lambda_2JE_i,\\
&h(E_1,E_\alpha)=\lambda_3JE_\alpha,~h(E_i,E_\alpha)=0.
\end{aligned}
\right.$$ Moreover, a Lagrangian immersion  $\psi:M\to
\mathbb{CH}^n(-4)$ satisfying the above conditions is  Hamiltonian
minimal if both $\psi_1$ and $\psi_2$ are Hamiltonian minimal, but
$\psi$ can never be minimal.
\end{thm}

\newpage

\end{document}